\documentclass[12pt]{article}
\usepackage{amsfonts, amsmath, amstext, amssymb}

\usepackage[dvips]{graphicx}
\usepackage{latexsym}

\newtheorem{thm}{Theorem}[section]
\newtheorem{prop}[thm]{Proposition}

\newtheorem{lemma}[thm]{Lemma}

\newtheorem{cor}[thm]{Corollary}

\newtheorem{definitiontemp}[thm]{Definition}
\newenvironment{defn}{\begin{definitiontemp}
\normalfont}{\end{definitiontemp}}

\def\bec{\begin{cor}}
\def\enc{\end{cor}}
\def\del{\delta}
\def\ep{\epsilon}

\def\bet{\begin{thm}}
\def\ent{\end{thm}}

\def\becor{\begin{cor}}
\def\encor{\end{cor}}

\def\bel{\begin{lem}}
\def\enl{\end{lem}}

\def\bedef{\begin{defn}}
\def\endef{\end{defn}}

\def\bep{\begin{prop}}
\def\enp{\end{prop}}

\newenvironment{pf}{\begin{trivlist}\item[\hskip\labelsep
{\it Proof.}]}{\end{trivlist}}

\newcommand{\set}[2]{\ensuremath{ \{ #1 : #2 \} }}

\renewcommand{\deg}[1]{\ensuremath{\text{deg}(#1)}}

\newcommand{\Z}{\mathbb{Z}}
\newcommand{\F}{\mathbb{F}}
\newcommand{\Q}{\mathbb{Q}}

\newcommand{\C}{\mathcal{C}}

\newcommand{\CC}{\mathfrak{C}}
\newcommand{\DD}{\mathfrak{D}}
\newcommand{\TFAB}{\mathfrak{Ab}_{TF}}
\newcommand{\FF}{\mathcal{ALG}}

\newcommand{\LL}{\mathfrak{L}}
\newcommand{\TT}{\mathcal{FBT}}

\def\SS{{\mathfrak S}}
\def\AA{{\mathcal{ACF}}}

\newcommand{\A}{\mathcal{A}}
\newcommand{\B}{\mathcal{B}}
\newcommand{\D}{\mathcal{D}}

\newcommand{\T}{\mathcal{T}}
\newcommand{\U}{\mathcal{U}}
\newcommand{\V}{\mathcal{V}}
\newcommand{\Vbar}{\overline{\mathcal{V}}}
\newcommand{\W}{\mathcal{W}}
\newcommand{\X}{\mathcal{X}}

\newcommand{\Qbar}{\overline{\mathbb{Q}}}

\newcommand{\Etilde}{\widetilde{E}}
\newcommand{\Ftilde}{\widetilde{F}}

\newcommand{\Ktilde}{\widetilde{K}}

\renewcommand{\L}{\mathcal{L}}

\newcommand{\G}{\mathcal{G}}

\newcommand{\I}{\mathcal{I}}
\newcommand{\avec}{\vec{a}}

\newcommand{\xvec}{\vec{x}}

\newcommand{\ACF}[1]{\ensuremath{\textbf{ACF}_{#1}}}

\newcommand{\comment}[1]{}
\newcommand{\Gal}[2]{\text{Gal}(#1/#2)}
\newcommand{\la}{\langle}
\newcommand{\ra}{\rangle}

\def\diverges{\!\uparrow}

\newcommand{\at}{\char'100}

\newcommand{\Eperm}{\ensuremath{E_{\text{perm}}}}

\newcommand{\Eset}{\ensuremath{E_{\text{set}}}}
\newcommand{\Ecard}{\ensuremath{E_{\text{card}}}}

\newcommand{\qed}{\hbox to 0pt{}\nobreak\hfill\rule{2mm}{2mm}}

\newcommand{\dom}[1]{\text{dom}(#1)}

\def\bfz{\boldsymbol{0}}
\def\s01{\ensuremath{\Sigma^0_1}}
\def\d02{\ensuremath{\Delta^0_2}}
\def\phi{\varphi}

\def\res{\!\!\upharpoonright\!}

\begin{document}

\title{Isomorphism and Classification\\for Countable Structures}
\author{Russell Miller
\thanks{This project grew out of the extended abstract \cite{MCiE},
although the results presented there and here are disjoint.
The author was supported by Grant \# DMS -- 1362206
from the National Science Foundation, and by several grants
from the PSC-CUNY Research Award Program and the Queens College
Research Enhancement Fund.
He wishes to acknowledge useful conversations with
Wesley Calvert, Rehana Patel, and Arno Pauly.}
}
\maketitle

\begin{abstract}
We introduce a topology on the space of all isomorphism types
represented in a given class of countable models, and use this
topology as an aid in classifying the isomorphism types.
This mixes ideas from effective descriptive set theory and computable
structure theory, extending concepts from the latter beyond computable
structures to examine the isomorphism problem on arbitrary countable
structures.  We give examples using specific classes of fields and of trees,
illustrating how the new concepts can yield classifications that reveal differences
between seemingly similar classes.  Finally, we use a computable homeomorphism
to define a measure on the space of isomorphism types of algebraic
fields, and examine the prevalence of relative computable categoricity
under this measure.
\end{abstract}

\section{Introduction}
\label{sec:intro}

Computable structure theorists have often restricted their study
to computable models of the particular theory or class in which
they are interested.  Attempts to measure the difficulty of determining
whether two given structures are isomorphic, such as
\cite{C04,GK02}, frequently consider only computable structures.
Indeed, sometimes the goal was to distinguish the computably presentable
structures from the others.
\cite{CCKM04}, \cite{CK06}, and \cite{KMV07} went further, defining computable
transformations of one class of countable structures into another, a notion
that we will adapt here to a broader context.  The question
of actually computing an isomorphism, in case the structures are isomorphic,
was initially posed for computable structures, as the notion of
computable categoricity, although its subsequent extension
to all countable structures (as \emph{relative computable categoricity})
was shown in \cite{AKMS89,Chis} to have a far simpler syntactic characterization
than computable categoricity itself has.  In \cite{DMN18, GK02, LMS18}.
questions of classification up to isomorphism are also restricted to
computable structures.

Of course, the class of \emph{all} structures isomorphic to a
given structure is a proper class, and hence unwieldy for study.
Here, however, given a class of countable structures, we will
consider all those representatives of the class which have the
domain $\omega$.  This is the natural context in which to study
questions of Turing computability on countable structures,
since $n$-ary functions and relations in the structure may then be
viewed as subsets of $\omega^n$, susceptible to analysis by
computability theorists.  One might request that subsets of $\omega$
be allowed as domains, but this essentially only matters when we
wish to include finite structures in our classes; otherwise, given the
domain itself, we may compute an isomorphism from the structure
onto a structure whose domain is all of $\omega$.

In this article we will address several classes of structures with domain
$\omega$.  First we examine $\FF_0$, the class of all algebraic field extensions of the field
$\Q$ of rational numbers.  From there we move on to algebraically closed
extensions of $\Q$ (of arbitrary transcendence degree $\leq\omega$),
which form a class usually viewed as simpler than the algebraic fields.
Last we consider the class $\TT$ of finite-branching trees, under the predecessor
function.  This class bears distinct
similarities to $\FF_0$, 
but in the context of all countable models, we will be able to draw distinctions
between them which were not so clear when only computable models were
considered.  The goal is to develop and illustrate the usefulness of studying
all countable models in the class.  We have deliberately chosen fairly simple
classes for this purpose, but nevertheless, several intriguing questions
arose out of these investigations, which are discussed in the final three sections.
It is hoped that subsequent work will extend this study
to classes with more complex isomorphism problems.

\section{Known Results on Computable Structures}
\label{sec:background}

Here we set up our notation and review some relevant background, mainly about computable structures.
We give all specific definitions for the class of algebraic field extensions of $\Q$;
those for algebraically closed fields and finite-branching trees are analogous.
Recall that the \emph{atomic diagram} of a structure $\A$ is the set of all atomic
sentences true in $\A$, in the language of $\A$ extended by new constants for elements
of the domain (which for us is always $\omega$).  For all structures on the domain
$\omega$ in this signature, we use a single effective bijective
G\"odel coding to view the atomic diagram as a subset of $\omega$.
\begin{defn}
\label{defn:notation}
\begin{align*}
&\FF_0 = \set{D\in 2^\omega}{D\text{~is the atomic diagram of a field algebraic over~}\Q}.\\
&\text{For $D\in\FF_0$,}\text{ $F_D$ is the field with domain $\omega$ and atomic diagram $D$.}\\
&\I(\FF_0) = \set{ (D_0,D_1)\in (\FF_0)^2}{F_{D_0}\cong F_{D_1}}.\\
&\FF_0^c = \set{e\in\omega}{(\exists D\in\FF_0)~\chi_D=\phi_e}.\\
&\text{For $e\in\FF_0^c$,}\text{ $F_e$ is the field $F_D$ for which~}\chi_D=\phi_e.\\
&\I(\FF_0^c) = \set{(e_0,e_1)\in (\FF_0^c)^2}{F_{e_0}\cong F_{e_1}}.\\
&\AA_0 = \set{D\in 2^\omega}{D\text{~is the atomic diagram of a model of }\ACF0}.\\
&\TT = \set{D}{D\text{~is the atomic diagram of an infinite finite-branching tree}}.
\end{align*}
\end{defn}
The language for the trees in $\TT$ is the language with equality and one unary
function symbol $P$, the \emph{predecessor function}.  The root $r$ of the tree
is the unique node satisfying $P(r)=r$, while for every other $x$, $P(x)$ is the
immediate predecessor of $x$ in the tree.  Of course, in this language, the greatest
possible height for a tree is $\omega$, and in fact all trees in $\TT$ have height $\omega$,
by K\"onig's Lemma, since they are infinite but finite-branching.
The tree is finite-branching if, for each node $x$, only finitely many $y$ satisfy $P(y)=x$.
Finite-branchingness for trees
cannot be stated in the first-order languages we have given; nor can algebraicity for fields.
The $0$ in $\FF_0$ and $\AA_0$ is to specify the characteristic.  We will not discuss
fields of positive characteristic here, but they behave in essentially the same way,
except that in $\FF_p$ it would be helpful to allow finite structures.

\comment{
For us a graph is really a symmetric (i.e., undirected) irreflexive graph, in the
language with equality and a binary edge relation $E$.  To get a \emph{pointed graph},
we add a single constant symbol to the language, to name a single node in the graph.
The graph is finite-valence if, for each node $x$, only finitely many $y$ satisfy
$E(x,y)$.
}

It is well-known that $\I(\FF_0^c)$ and the analogously defined set $\I(\TT^c)$ 
are both $\Pi^0_2$-complete.  The fact that they are $\Pi^0_2$ stems from the following
characterizations, which essentially are consequences of K\"onig's Lemma.
\begin{lemma}
\label{lemma:Pi2}
Two algebraic fields of equal characteristic are isomorphic if and only if every
finitely generated subfield of each one embeds into the other (equivalently, if
every polynomial over the prime subfield with a root in either field
also has a root in the other).  Two finite-branching trees are isomorphic if and
only if every finite subtree of each one embeds into the other.  
\end{lemma}
Moreover, these equivalence relations (on the relevant subsets
of $\omega$, namely $\FF_0^c$ and $\TT^c$) are computably reducible
to each other, and also computably bireducible with the equivalence relation
$=^e$ on $\omega$ defined by
$$ i =^e j \iff W_i = W_j.$$
This means that there are computable functions $f$ (which is total, with image
$\subseteq\FF_0^c$) and $g$ (with domain $\supseteq\FF_0^c$) such that
\begin{align*} (\forall i,j)~[i =^e j&\iff F_{f(i)}\cong F_{f(j)}]\\
\&~~~(\forall i,j\in\FF_0^c)~[F_i\cong F_j&\iff g(i) =^e g(j)],
\end{align*}
and likewise for $\TT^c$.
These functions $f$ and $g$ are \emph{reductions} between the equivalence
relations involved, in the sense described in \cite{MSEALS}, for example.
The function $f$ here is defined by setting $f(i)$ to be the index for an atomic diagram
of the field $\Q[\sqrt{p_n}~:~n\in W_i]$, where $p_n$ is the $n$-th prime number.
To define the function $g$, fix a listing $h_0,h_1,\ldots$ of all irreducible
polynomials in $\Q[X]$, and let $W_{g(i)}=\set{n\in\omega}{h_n\text{~has a root in~}F_i}$.
Likewise, to reduce $=^e$ to $\I(\TT^c)$, let $f(i)$ be the index of the atomic diagram
of the tree with a single infinite ``spine'' of length $\omega$ and with one other node
at level $n+1$ (above the spinal node at level $n$) just if $n\in W_i$; for the opposite
reduction, use a computable list of all finite trees and consider the set of those which embed
into the given tree.  

\section{Algebraic Fields}
\label{sec:fields}

The reductions described above between $\FF_0^c$ and $\TT^c$ were
computable reductions on equivalence relations on subsets of $\omega$.
We now switch to the notion of a computable reduction from one equivalence
relation on one subset of Cantor space to another equivalence relation on
another such subset, as has been done in \cite{MSEALS} and other places.
This can be seen as a generalization of \cite[Defn.\ 8]{CK06}, where it is done
for the relation of isomorphism on classes such as $\FF_0$ and $\TT$.
Here we will allow ourselves to work with arbitrary equivalence relations
on Cantor space or Baire space, and our goal will be not merely reductions,
but actual homeomorphisms.  To begin with, therefore, we explain the topologies
in question, with $\FF_0$ as our first example.

The sets $\FF_0$, $\AA_0$, and $\TT$ are all subsets of Cantor space,
so each one becomes a topological space in its own right, under the subspace topology.
The basic open sets in the topology on $\FF_0$, for instance, are defined to be the sets
of the form
$$ \U_{\sigma}=\set{D\in\FF_0}{\sigma\subset D},$$
with $\sigma$ ranging over all initial segments of elements of $\FF_0$ (or over all of $2^{<\omega}$,
in which case many of the sets above are empty).  Then, under the equivalence relation
of isomorphism on $\FF_0$, defined by $\I(\FF_0)$, we get the quotient topology on
the quotient space $\FF_0/\!\cong$, sometimes known as the \emph{isomorphism space}
of $\FF_0$.  By definition this means that a subset $\V$ of $\FF_0/\!\cong$ is open there
if and only if its preimage under the quotient map
$$ \U = \set{D\in\FF_0}{\text{the isomorphism type of $F_D$ lies in~}\V}$$
is open in the topology on $\FF_0$.  Since such an open set $\U$ must be the
union of certain basic open sets $\U_\sigma$, we see that
an open set $\V$ in $\FF_0/\!\cong$
is defined by taking any set of initial segments $\sigma$ of various $D\in\FF_0$ and
letting $\V$ contain all isomorphism types of fields $F_D$ with $D\in\U_\sigma$,
for all these strings $\sigma$.  That is, $\V$ contains all isomorphism types
such that one of these $\sigma$ can be an initial segment of a field of that type.
(The analogous statement holds for any other class of structures on the domain $\omega$, of course.)
%

Now a string $\sigma$ may contain
enough information to ensure that all its extensions in $\FF_0$ contain roots of one
or of several polynomials:  it might state, for example, that the square of some element,
when added to the multiplicative identity, equals the additive identity, i.e., that $X^2+1$ has a root.
It might also give incomplete information.  For example, $\sigma$ could
specify that the field contains elements $x$ and $y$ with $x^8=y^{12}=2$,
but without specifying the minimal polynomial of $y$ over $\Q(x)$.
Since $Y^{12}-2$ factors as $(Y^3-x^2)(Y^3+x^2)(Y^6+x^4)$ over $\Q(x)$,
we cannot determine even the degree of the extension dictated by $\sigma$,
let alone its isomorphism type.  However, there are only finitely many possibilities,
and in any particular field, some finite piece
of the atomic diagram must specify the minimal polynomial of $y$ over $\Q(x)$,
so in practice we normally just choose a string sufficiently long to specify
the isomorphism type of the relevant subfield.

Since a single $\sigma$ is finite, it can only specify finitely many extensions of finite
degree to be contained in the field.  So an open set $\V$ in $\FF_0/\!\cong$
will contain every isomorphism type realizing the finite-degree extensions
dictated by $\sigma$, as $\sigma$ ranges over some set of strings.
Of course, by the Primitive Element Theorem, every finite-degree algebraic extension
of $\Q$ is generated by a single element, so it is natural to think of $\V$ as 
being defined by a particular subset $S$ of $\Q[X]$, corresponding
to the set of strings $\sigma$ defining its preimage $\U$:  $\V$ is the set
of all isomorphism types containing a root of any polynomial in $S$.
(For structures more generally, $\V$ will contain every isomorphism type realizing
any of the finite partial quantifier-free $n$-types defined by these strings $\sigma$.)
Conversely, for any given finite-degree extension
$K$ of $\Q$, we can choose a primitive generator $y$ of $K$ and its minimal
polynomial $p(Y)$ over $\Q$.  Take the open subset of $\FF_0$ given by
$$ \U=\bigcup_\sigma \U_\sigma,$$
where the union is over all $\sigma$ which state that, for some $y<|\sigma|$, $p(y)=0$.
Thus the image $\V$ of $\U$ under the quotient map is an open set (since $\U$
is closed under isomorphism) and will contain exactly those algebraic fields extending $K$.
This yields the following.
\begin{lemma}
\label{lemma:FF_0}
The topology on $\FF_0/\!\cong$ has a basis consisting of those sets defined by
the existence of a root for one irreducible polynomial from $\Q[Y]$.
Equivalently, these are the sets of the form $\set{[D]\in\FF_0/\!\cong}{K\subseteq F_D}$,
as $K$ ranges over all finite-degree field extensions of $\Q$.
\end{lemma}
It is important to notice that these sets do form a basis, not just a subbasis.
To understand the potential difficulty, consider the sets $\V$
and $\W$ containing, respectively, those fields with an eighth root of $2$ and
those with a twelfth root of $2$.  As seen earlier, the intersection $\V\cap\W$,
which is nonempty, cannot be defined by the existence of any single subfield;
there are three fields $K_0$, $K_1$, and $K_2$, pairwise incomparable under
embedding, such that the fields in $\V\cap\W$ are those containing at least one of these three
as a subfield, but it is impossible to give two fields, let alone just one, which
characterize $\V\cap\W$ in the same way.  (The three fields correspond to
the three factors of $Y^{12}-2$ over $\Q(x)$, seen above.).  However,
if $\V_i$ contains all isomorphism types of fields extending $K_i$, then each
$\V_i$ is an element of the basis in Lemma \ref{lemma:FF_0}, and $\V\cap\W=\V_0\cup\V_1\cup\V_2$.


\begin{cor}
\label{cor:FF_0}
In the topology on $\FF_0/\!\cong$, the only open set containing the isomorphism type of $\Q$
is the entire space $\FF_0/\!\cong$, whereas the type of the algebraic closure $\Qbar$
lies in every nonempty open set.
\qed\end{cor}
This corollary is not unique to $\FF_0$; it would hold in any class $\CC\subseteq 2^\omega$
of structures, provided that $\CC$ contains both a least and a greatest
structure under embedding.  It can readily be related to \cite[Prop.\ 4.1]{CK06},
re-expressed in terms of topology.
As a topological space, $\FF_0/\!\cong$ is not merely compact:
every open cover has a subcover of size $1$, since some element
of the cover is an open set containing the type of $\Q$.
However, the corollary does establish that this topology
on $\FF_0/\!\cong$ is not homeomorphic to any of the best-known topological spaces
of size continuum.  We discuss this further in Section \ref{sec:Scott topology}.

To address this shortcoming, we now adjoin to the language $\L$ of fields, for each $n>1$,
the $n$-ary \emph{root predicate} $R_n$, defined as follows:
$$ \models_{\F} R_n(a_0,\ldots,a_{n-1}) \iff (\exists x\in\F) x^n+a_{n-1}x^{n-1}+\cdots+a_0=0.$$
\begin{thm}
\label{thm:algfields}
In the language $\L^*$ with the predicates $R_n$, the class $\FF_0^*/\!\cong$
of all algebraic fields of characteristic $0$, modulo isomorphism, is computably
homeomorphic to Cantor space $2^\omega$ (under its usual topology).
In particular, it is compact.
\end{thm}
In this context, a homeomorphism $h$ from $\FF^*_0/\!\cong$ onto $2^\omega$
is \emph{computable} if there is one Turing functional $\Gamma$ such that, for every
$D\in\FF^*_0$, $\Gamma^D=h([D])$, the image of the equivalence class of $D$,
and a second Turing functional $\Phi$ such that, for every $D\in\FF^*_0$,
$\Phi^{h([D])}$ lies in the class $[D]$.  Every homeomorphism is
continuous, hence $S$-computable for some $S\subseteq\omega$, but we assert here
that for this homeomorphism, no oracle is required.  (We refer readers to \cite{W00}
for details of the relation between computability and continuity.)

\begin{pf}
Fix a computable copy $E$ of the algebraically closed field $\Qbar$,
and let $f_0,f_1,\ldots$ enumerate all monic polynomials in $E[X]$ of degree
$>1$.  For each $\sigma\in 2^{<\omega}$, we now define a subfield $F_\sigma\subseteq E$
and a polynomial $f_\sigma$ from this list.  The root node $\lambda$ of $2^{<\omega}$
has $F_\lambda=\Q$, the prime subfield of $E$.  For each $\sigma$, starting with $\lambda$, 
we now define $f_\sigma$ recursively to be $f_n$ for the least $n$ such that:
\begin{itemize}
\item
$f_n$ lies in $F_\sigma[X]$ (within $E[X]$) and is irreducible of prime degree there; and
\item
for the least root $x$ of $f_n$ in $E$, and for every $\tau$
with $\tau\widehat{~}0\subseteq\sigma$, the subfield $F_\sigma(x)$ of $E$
contains no root of $f_\tau$.
\end{itemize}
A theorem of Kronecker from \cite{K1882} gives us splitting algorithms for every finitely
generated subfield of $E$, uniformly in the generators, so these conditions
are both decidable.  (For details see \cite[Thm.\ 2.1]{M10}.)
Since $F_\sigma$ is finitely generated, there must exist
some $f_n$ satisfying these conditions.  This defines $f_\sigma$,
and we now set $F_{\sigma\hat{~}0}=F_\sigma$ and
$F_{\sigma\hat{~}1}=F_\sigma(x)$, for the least root $x$ of $f_\sigma$ in $E$.
(The roots of $f_{\sigma}$ are elements of the domain $\omega$ of $E$,
so by the ``least root,'' we mean the smallest element of $\omega$ which is a root.)
This completes the recursion, defining $F_\sigma$ and $f_\sigma$ for every $\sigma$.

We now define the functionals $\Gamma$ and $\Phi$, which will be inverses
of each other when viewed as maps between $\FF^*_0/\!\cong$ and $2^\omega$.
$\Gamma$ accepts as input the atomic diagram $D$ of any field $F_D$ in $\FF^*_0$,
in the language with the predicates $R_n$.
Starting with the empty string $\sigma_0$ and the unique embedding $g_0$ of $F_{\sigma_0}=\Q$
into $F_D$, it defines each string $\sigma_{m+1}\in 2^{m+1}$ to extend $\sigma_m$ by setting:
\begin{itemize}
\item
$\sigma_{m+1}(m)=1$ if $F_D$ contains a root of $f_{\sigma_m}$, in which case
$g_{m+1}$ extends $g_m$ by mapping the least root of $f_{\sigma_m}$ in $F_{\sigma_{m+1}}$
to the least root in $F_D$ of the $g_m$-image of $f_{\sigma_m}$ in $F_D[X]$;
\item
$\sigma_{m+1}(m)=0$ if not, in which case $g_{m+1}=g_m$.
\end{itemize}
So in both cases $g_{m+1}$ embeds $F_{\sigma_{m+1}}$ into $F_D$.
With the root predicate in the language, this procedure is effective:  $D$ allows $\Gamma^D$ to
determine the correct value for $\sigma_{m+1}(m)$, extending $\sigma_m$
to length $m+1$.  The output of $\Gamma$ is the function $h=\cup_m \sigma_m$,
and $g=\cup_m g_m$ is an isomorphism from the subfield $\cup_m F_{h~\!\res~\!m}$ of $E$
onto $F_D$.
Notice that for each $m$, the $g$-image of $f_{h~\!\res~\!m}$ in $F_D[X]$ has a root
in $F_D$ if and only if $h(m)=1$.

$\Phi$ accepts any $h\in 2^\omega$ as oracle, and defines a field $K=\cup_m F_{h\upharpoonright m}$,
using the fields $F_\sigma$ defined above.  This is done in such a way that $\omega$
is the domain of $K$ (even if $h$ is the constant function $0$.)  Moreover, with $h$ as oracle,
$\Phi^h$ will eventually decide whether a given polynomial $X^d+a_{d-1}X^{d-1}+\cdots+a_0$
in $K[X]$ has a solution in $K$, because this polynomial (or factors of it) will eventually
appear as polynomial(s) $f_{h\upharpoonright m}$ for various $m$, and the existence of a root of
$f_{h\upharpoonright m}$ in $K$ will be determined by the value of $h(m)$.
Thus $\Phi^h$ decides the atomic diagram of the field $K$ it builds, even in the expanded
language with the root predicates.  Moreover, applying $\Gamma$ to the atomic diagram
of this field will clearly yield the same $h$ we started with, and likewise,
applying $\Phi$ to the function $\Gamma^D$ will yield a field isomorphic
to $F_D$.  So the two computable functions are inverses, and we have
the desired homeomorphism.
\qed\end{pf}

Notice that this proof would not work if we did not have the root predicates
in our language.  As noted above, the structure $\Qbar$
lies in every nonempty open set of $\FF_0/\!\cong$, so clearly $\FF_0/\!\cong$
and $\FF_0^*/\!\cong$ cannot be homeomorphic.
Having the root relations in the language allows the atomic diagram of $\Q$
to differentiate itself from those of other algebraic fields, as the failure of
these relations to hold (on specific elements of $\Q$) expresses the absence
from $\Q$ of roots of certain polynomials.  There does exist a continuous
bijection from $\FF_0^*/\!\cong$ onto $\FF_0/\!\cong$, sending each isomorphism
type to itself; indeed, this bijection is computable, mapping the atomic diagram
of any $F_D^*$ to the atomic diagram of its reduct.  (It follows that $\FF_0/\!\cong$
is also compact, as we already knew.)  However, the inverse map is not continuous.

It is important to notice that, in order to make $\FF_0$ homeomorphic to
$2^\omega$, we had to be careful not to add too \emph{much} information
to the language either.  It would have been very natural to extend each
field $F$ in $\FF_0$ to its \emph{jump} $F'$, as defined in \cite{M09,SS09}
by Montalb\'an, Soskov, and Soskova.  The root predicates $R_n$
all belong to the signature of the jump of a field, and so we would still
have a computable function $\Gamma$ from the class $\FF_0^{(1)}/\!\cong$ of jumps
of fields in $\FF_0$, modulo isomorphism, onto Cantor space;
it would be exactly the same map as from $\FF_0^*$ onto $2^\omega$.
However, in $\FF^{(1)}_0/\!\cong$, the isomorphism type of $\Q$ forms an open
set all by itself, because the language of the jump of a field contains
a predicate for the following computable existential $L_{\omega_1\omega}$
statement about the field:
$$ (\exists p\in\Q[X])(\exists x)[p\text{~is irreducible of degree}>1~\&~p(x)=0].$$
A field $F\in\FF_0$ is isomorphic to $\Q$
just if the predicate for
this statement is false in $F'$.  Thus the topological space $\FF^{(1)}_0/\!\cong$
contains a singleton open set, whereas $2^\omega$ does not, and so
there can be no homeomorphism between them.  In particular,
the inverse of the map $\Gamma$ is not computable, nor even continuous.
As an insight into this failure, notice that, with the constant function $0$
as $h$, an inverse function $\Phi^h$ would have to compute a copy of the jump
of the field $\Q$, and this copy would have to compute the degree $\bfz'$, as all jumps
of structures do; but it is clearly impossible for a computable functional $\Phi$,
with the computable oracle $h$, to output an atomic diagram which decides $\bfz'$.
Of course, a continuous inverse function might require an oracle $S$,
i.e., might only be given by an $S$-computable functional $\Phi^S$,
but a similar argument with a field whose
spectrum is the upper cone above $\deg{S}$ shows this also to be impossible.

\begin{prop}
\label{prop:Scott}
Let $\CC\subseteq 2^\omega$ be a class of structures, closed under isomorphism.
Then, for each structure $\D$ with atomic diagram $D\in\CC$,
the singleton of the isomorphism type of $\D$ forms an open set
in $\CC/\!\cong$ if and only if there is some (finitary) $\Sigma_1$
formula $\exists\xvec\alpha$ in the language which is a Scott sentence
for $\D$ within the class $\CC$.  (This means that, among structures in $\CC$,
$\alpha$ holds just in those structures isomorphic to $\D$.)
\end{prop}
To clarify:  in the example with $\Q$ in $\FF_0^{(1)}$, the Scott sentence is actually a quantifier-free
sentence $\alpha$ in the language of jumps of fields.  In the more general context, however,
finitary $\Sigma_1$ formulas may be necessary.
\begin{pf}
First, if there is such a Scott sentence $\exists\xvec\alpha(x_1,\ldots,x_n)$ in the language,
then for every $\C\cong\D$ with $\Delta(\C)\in\CC$, we can
choose $(k_1,\ldots,k_n)\in\omega^n$ realizing $\alpha$ in $\C$
and let $\sigma_{\C}$ be an initial segment of $\Delta(\C)$ sufficiently long to state that
$\models_{\C}\alpha(k_1,\ldots,k_n)$.  Then $\U_{\sigma_{\C}}\cap\CC$ is an open
set of $\CC$ containing only diagrams of structures in which $\exists\xvec\alpha$ holds,
including $\D$ itself, and the union of all these $\U_{\sigma_{\C}}$ (over all $\C\cong\D$
in $\CC$) is an open set whose image in $\CC/\!\cong$ is the singleton of the class of $\D$.
Thus this singleton is an open set in $\CC/\!\cong$.

Conversely, suppose that the isomorphism type of $\D$ forms a singleton open set in
$\CC/\!\cong$.  Then its preimage $\U$ in $\CC$ is open and contains the atomic diagram $D$ of $\D$,
but contains no atomic diagrams of any structures $\not\cong\D$.  Let $\sigma\subset D$ be sufficiently
long that $\U_\sigma\cap\CC$ is contained in this open preimage, and let $\alpha$
be the conjunction of all atomic statements in $\sigma^{-1}(1)$ and all negations of
atomic statements in $\sigma^{-1}(0)$.  Replacing each element $j$ of $\omega$
used in $\alpha$ by the variable $x_j$ yields a finitary $\Sigma_1$ formula
$$ \exists x_0\cdots\exists x_{|\sigma|}~\alpha(x_0,\ldots,x_{|\sigma|})$$
which holds in $\D$.  Whenever $\C$ is a structure (with diagram in $\CC$)
realizing this $\Sigma_1$ formula, there must be some $\B\cong\C$ such that the
atomic diagram of $\B$ begins with $\sigma$:  just use a permutation of $\omega$
mapping the elements of $\C$ realizing $\alpha$ onto $(0,\ldots,|\sigma|)$.  Since
$\CC$ is closed under isomorphism, the diagram of $\B$
lies in $\U_\sigma\cap\CC$, hence lies in $\U$, forcing $\B\cong\D$.
Thus $\exists\xvec\alpha$ is indeed a Scott sentence for $\D$ within $\CC$.
\qed\end{pf}
We will discuss singleton open sets further in the next section.

\section{Algebraically Closed Fields}
\label{sec:ACF}

We postponed consideration of algebraically closed fields until after
algebraic fields because, although the former is a simpler class, the latter
is more representative of the usual situation.  An algebraically closed field
is determined up to isomorphism by its characteristic and its transcendence
degree; the field is countable just if the latter is $\leq\omega$.  So there are
only countably many isomorphism classes of countable algebraically closed fields.
We will restrict our discussion to those of characteristic $0$.

With only the usual language of fields, the goal is to count the number
of elements in a transcendence basis of an ACF $F$.  Recognizing transcendentals
is not a decidable process, so we first consider the easier situation in which
the language is enriched by $n$-ary \emph{dependence predicates} $D_n$,
for every $n>0$, which hold of an $n$-tuple from $F$ if and only if
that tuple is algebraically dependent over the prime subfield $\Q$.
Now a field $F$ (with domain $\omega$) has a natural transcendence
basis, containing exactly those elements $x$ which are transcendental over
all preceding elements of the domain $\omega$ of $F$:
$$ (\forall p\in\Q[Y_0,\ldots,Y_x])~[p(0,1,\ldots,x-1,x)=0\implies p(0,1,\ldots,x-1,Y_x)=0.$$
This basis is not invariant under automorphisms, of course, but
its elements are recognizable from the atomic diagram
(in the language with the predicates $D_n$).  So it is readily seen that
for the class $\AA_0^D$ of algebraically closed fields of characteristic $0$,
with domain $\omega$, in the language with the dependence predicates,
the space $\AA_0^D/\!\cong$ is homeomorphic to $2^\omega/\Ecard$, where
the equivalence relation $\Ecard$ (the \emph{cardinality relation}) is defined by:
$$ A~\Ecard~B \iff |A|=|B|.$$
In one direction, a field $F$ is mapped to the transcendence basis defined above;
while in the reverse direction, each time we find a new element of the oracle set $A$,
we adjoin a new transcendental to the field we are building, while always taking
one further step to make the field algebraically closed.  Notice that this does allow
computation of a transcendence basis of the field (namely, the set of elements
added in response to the discovery of new elements of $A$), from which in turn
one can compute the dependence relations on the field.

The cardinality relation is not among the usual Borel equivalence relations on $2^\omega$.
It has only countably many equivalence classes, among which one is a singleton
(the class of $\emptyset$), one has size continuum,
and all others are countable.  $\Ecard$ is not \emph{smooth}:
the equality relation on $2^\omega$ has no Borel reduction
to $\Ecard$, as $\Ecard$ has too few classes.  In the topology on
$2^\omega/\Ecard$, the basic open sets are those defined by having at least $n$ elements,
and these are the only nonempty open sets in the space:
$$ \U_n = \set{A\in 2^\omega}{|A|\geq n}.$$
Therefore, the $\Ecard$-class of infinite sets belongs to every nonempty open set,
whereas the $\Ecard$-class of $\emptyset$ belongs to no open set except the entire space.
The same holds in $\AA_0^D/\!\cong$:  the prime model $\Qbar$
belongs to no open set except the whole space, whereas the algebraically closed
field with transcendence degree $\omega$ belongs to every nonempty open set.  
The open sets are just those closed upwards under transcendence degree.

Of course, by adding more definable relation symbols to the language, one could
change this topology.  The transcendence degree itself is $L_{\omega_1\omega}$-definable,
after all.  One specific method is to add $n$-ary predicates $B_n$, which hold
of an $n$-tuple $\avec$ if and only if $\avec$ is a transcendence basis for $F$.
(It is helpful to include $B_0$, which states $\forall x D_1(x)$; this is essentially $0$-ary,
and holds only in transcendence degree $0$.)  In the space $\AA_0^*/\!\cong$
of algebraically closed fields in this language, the isomorphism
type of a field of finite transcendence degree will itself form an open set, but the space
is not quite discrete, as the only open sets containing the type of infinite transcendence
degree are those sets which, for some $d$, contain all types of degree $\geq d$.
Thus this space is still compact.  To make the singleton of the infinite-degree type
open as well (and thus make the space noncompact), one must add a predicate
$C$ defining it, such as $(\forall d)(\exists\xvec)\neg D_d(\xvec)$.
Once again this is only $L_{\omega_1\omega}$-definable.

Several principles are illustrated here.  First, this process recalls Proposition
\ref{prop:Scott}:  in a space $\DD/\!\cong$, the singleton set
containing just the isomorphism type of a structure $\A$ from $\DD$ is an open set
if and only if $\A$ has a Scott sentence which is (finitary) $\Sigma_1$ in the language of $\DD$.
Such a sentence says that among structures in $\DD$, being isomorphic to $\A$ is equivalent
to containing a finite collection of elements satisfying a certain configuration, and so
a finite piece of the atomic diagram of $\A$ is sufficient to guarantee isomorphism to $\A$.
Each predicate $B_n$ above yielded a Scott sentence for the algebraically closed field
of transcendence degree $n$ (among algebraically closed fields of characteristic $q$),
but the infinite-degree field required a bit more work.

Second, once all the Scott sentences (including for infinite degree) were adjoined to the
language, the space became the discrete countable topological space, hence noncompact.
Consider the general situation of a class $\DD$ of countable structures and a class $\DD^*$
of the same structures in a language extended by predicates
which are $L_{\omega_1\omega}$-definable in the original language.
Thus $\DD/\!\cong$ and $\DD^*/\!\cong$ have an obvious bijective correspondence as sets,
and the map from $\DD^*/\!\cong$ onto $\DD/\!\cong$ is clearly computable, since the atomic
diagram of a reduct can easily be computed from the atomic diagram of
the structure.  (If the languages were countable but noncomputable,
then this reduction might not be computable either, yet it is still continuous, requiring
only an oracle for the languages.)  However, the inverse map, from the reduct $\A$ to
the original $\A^*$, may well not be computable, and often is not even continuous,
as seen in the examples above, where the topological spaces $\DD^*/\!\cong$ and $\DD/\!\cong$
are not homeomorphic.  Therefore, although compactness of the space in the larger language
implies compactness of the space in the smaller language (as the continuous image of a
compact space is always compact), the converse fails.  The class $\AA_0$ and
its expansions provide a ready counterexample.

The situation above resembles the well-known case of the \emph{jump map}, taking a set
$A\in 2^\omega$ to its jump $A'$.  This map is not continuous -- its points of continuity
are precisely the $1$-generic sets -- but its inverse is readily seen to be a continuous map
on the domain $\set{A'}{A\in 2^\omega}$, and in fact extends to a continuous map on
all of $2^\omega$.

\section{Finite-Branching Trees and Graphs}
\label{trees}

For us, the language of trees has equality and a single unary function symbol $P$,
the \emph{predecessor function}.  We define $P(r)=r$ for the root $r$ of the tree,
for convenience and also to enable the root to be identified from the atomic diagram.
We stipulate that for every $x$ in $\T$, there must exist
an $n$ with $P^{n+1}(x)=P^n(x)$, meaning that every node is only finitely many
levels above the root.  (The \emph{level} of $x$ is the least $n$ for which this holds.)
A tree $\T$ is \emph{finite-branching} if, for every $x\in\T$, the inverse image
$P^{-1}(x)$ is a finite set (possibly empty).  These definitions use $L_{\omega_1\omega}$
formulas, of course; they cannot be expressed finitarily.  The class of all infinite
finite-branching trees with domain $\omega$ is denoted by $\TT$.

We enrich the class to $\TT^*$ by adjoining to this language, for each $n>0$,
the unary \emph{branching predicate} $B_n$,
defined as follows:
$$ \models_{\T} B_n(a) \iff (\exists
x_0,\ldots,x_{n-1}\in\T)\left[\bigwedge_{i<n} P(x_i)=a
~\&\bigwedge_{j<k<n} x_j\neq x_k\right].$$
Since these trees are all finite-branching, it would be equivalent
to have predicates saying that $a$ has exactly $n$ immediate successors.
\begin{thm}
\label{thm:FBtrees}
The isomorphism space $\TT^*/\!\cong$ of infinite finite-branching trees,
in the language with these predicates $B_n$, is
homeomorphic to Baire space $\omega^\omega$.
Thus it is not compact and not homeomorphic to $\FF_0^*/\!\cong$.
\end{thm}
\begin{pf}
We write $T|_m$ for the subtree of $T$ containing all nodes at levels $\leq m$.
Fix a computable list $\L=\{ S_0,S_1,\ldots\}$ of all finite trees,
without repetitions.  Given the atomic diagram of a tree
$T$ in $\TT$, the predicates $R_n$ allow us to determine, for each $m$,
exactly which tree $S_{f(m)}$ in $\L$ is isomorphic to the finite subtree $T|_m$.
(Without the branching predicates, we would not be able
to determine this.)  So our functional $\Gamma$ outputs a function $h=\Gamma^{\Delta(T)}$
in $\omega^\omega$, defined by setting $h(m)$ to be the size of the set
$$ \set{j<f(m+1)}{S_j\text{~is a tree of height $m+1$ with $S_j|_m\cong T|_m$}}.$$
The intuition is that we list out and number all trees in $\L$ of height
exactly $m+1$ which extend $T|_m$ to the next level.
These are all the possibilities for $T_{m+1}$.  We define $h(m)$ so that
the tree isomorphic to $T|_{m+1}$ is the $h(m)$-th tree on the list (allowing $h(m)=0$).

The inverse $\Phi$ is also defined readily:  given $h\in\omega^\omega$,
$\Phi^h$ builds a tree $T$ by choosing the $h(0)$-th possibility among all trees
on $\L$ of height $1$, then extending it to the $h(1)$-st possibility among all trees
of height $2$ on $\L$ which extend the tree of height $1$ already built, and so on.
Since $\Phi^h$ knows from this the exact isomorphism type of the first $m$ levels of $T$,
it can compute the branching predicates for this tree and include them in the atomic diagram.
Clearly this procedure is inverse to the $\Gamma$ above.
\qed\end{pf}

The main point of Theorem \ref{thm:FBtrees} is the remark at the end:
that $\TT^*$ and $\FF_0^*$, modulo isomorphism, are not homeomorphic,
and hence do not have the same classification.  In light of the substantial similarities
between the classes $\FF_0^c$ and $\TT^c$ of computable structures (including
the analogous use of the root and branching predicates) this could be a surprise.
Considering all countable structures sheds light on the situation in a way that
we did not get from the computable structures alone.

The situation becomes a little more complicated when we extend $\TT$ to allow
finite trees:  now it is possible that a tree $T$ simply ends at some finite level. 
The simplest device for allowing finite structures is to allow the atomic diagram
to include or exclude the statement $n=n$ for each $n\in\omega$; the domain of
the structure, which now consists of those $n$ satisfying $n=n$, is thus decidable from
the atomic diagram.  Of course, the finiteness of the structure is not decidable.
With branching predicates, however, finiteness is $\Sigma_1$, since the branching predicates
will tell us when we have reached the top level of a finite tree $T$.  Without giving all the
details, we state here that this class of structures, under isomorphism, is homeomorphic
to the topological space which includes all of $\omega^{\omega}$ and also includes,
for each $\sigma\in\omega^{<\omega}$, one path $h_\sigma$ which extends $\sigma$
but is isolated above $\sigma$.  (Equivalently, it is homeomorphic to the subclass $\set{h\in\omega^\omega}{
(\forall n)[h(n)=0\implies h(n+1)=0]}$, under the subspace topology inherited from $\omega^\omega$.)

The class $\TT$ is closely related to the class of finite-valence infinite connected pointed graphs,
i.e., connected symmetric irreflexive graphs for which each node has only finitely many neighbors,
with one extra constant in the language (which makes the graph \emph{pointed}).
Indeed, this class turns out also to be classified effectively via a computable homeomorphism
onto Baire space, once the valence function is added to the language.
The \emph{valence function} tells how many neighbors each node has;
equivalently, one can add unary predicates $V_n(x)$ saying that $x$ has
$n$ neighbors.  The homeomorphism is similar to that for $\TT^*$:
$f(0)$ describes the isomorphism type of the constant node and its neighbors;
$f(1)$ describes the isomorphism type of their neighbors, i.e., of the subgraph of radius $2$
around the constant, given the subgraph of radius $1$ described by $f(0)$; and so on.

When one abandons the constant, leaving the class of connected finite-valence graphs,
the isomorphism problem becomes $\Sigma^0_3$ instead of $\Pi^0_2$.
This is a class of structures for possible further study.  In general, as the isomorphism
problem grows more difficult, we suspect that it will be necessary to use various
of the well-known basic equivalence relations on $2^\omega$ or on $\omega^\omega$,
such as the relation $E_0$ of having finite symmetric difference, or the relation
$\Eperm$ of having the same columns (when viewed as subsets of $\omega^2$)
up to permutation.  This means that a homeomorphism onto the space
$2^\omega/E_0$ or $\omega^\omega/\Eperm$, for instance, would fall under the head
of an effective classification, as would $2^\omega/\Ecard$, which we already used
to classify $\AA_0^D$.  Certain other useful equivalence relations will be
discussed in Section \ref{sec:Scott topology}.

\comment{
Additional equivalence relations on Cantor space
which seem likely to be useful include:
\begin{align*}
A~\Ecard~B &\iff |A|=|B|.\\
A =^e B &\iff \set{x}{\exists y~\la x,y\ra\in A} = \set{x}{\exists y~\la x,y\ra\in B}.  \\
A =^f B &\iff  (\forall x)~\set{y}{\la x,y\ra\in A}~\Ecard~\set{y}{\la x,y\ra\in B}.
\end{align*}
The latter two are both $\Pi^0_2$ relations, Borel-equivalent to equality on $2^\omega$
(and therefore not much studied in descriptive set theory) but not computably equivalent to it.
We can classify the subrings of $\Q$ effectively
via a homeomorphism onto $2^\omega/\!=^e$, while $=^f$ has applications
in classifying the class of all (finite or countable) equivalence structures (consisting
of a single equivalence relation on a domain $\subseteq\omega$.)
The relation $\Ecard$ is unusual in this context, having only countably many equivalence classes,
almost all of which are themselves countable, but we saw its usefulness in Section \ref{sec:ACF}.
}

\section{Measure and Category for Field Properties}
\label{sec:measure}

In classes of structures with a sufficiently nice classification, we can consider
various properties of isomorphism types and inquire into the frequency with which
those properties occur.  For an example, we return now to the class $\FF_0^*$
of algebraic fields of characteristic $0$, in the language with all root predicates $R_n$.
By Theorem \ref{thm:algfields}, we have a computable homeomorphism $\Phi$
from Cantor space $2^\omega$ onto the space $\FF^*_0/\!\cong$ of all
isomorphism types of such fields, allowing the transfer of the standard
notions of Lebesgue measure and Baire category from $2^\omega$ onto $\FF^*_0/\!\cong$.
The Lebesgue measure of a class of isomorphism types of fields is simply the
Lebesgue measure of the corresponding set of reals in Cantor space,
and the class of fields is nowhere dense, meager, etc., just if the corresponding
set of reals is.
Later, in Section \ref{sec:Haar}, we will discuss whether Lebesgue measure is the best measure
to use for this purpose.  First, however, we illustrate our goals by offering potential uses
of measure and category.

It is clear that the class of all isomorphism types of fields normal over $\Q$ has measure $0$,
and is meager.  Infinitely many irreducible polynomials $f$ in $\Q[X]$ have non-cyclic
Galois groups, and for each of these, there is positive probability that $f$ will have
at least one root but not all of its roots in a given field.  (More precisely, one can give infinitely
many polynomials whose Galois groups are all $S_3$, the symmetric group on the three roots
in $\Qbar$, and such that the splitting field of each is linearly disjoint from that of the others.
Each of these polynomials has the same positive probability of having exactly one root
in a given field, and linear disjointness ensures that these probabilities are all
independent:  the number of roots of one such polynomial in $F$ is independent of the number of
roots of any of the others.  We will not go further into the details here; the reader may refer to
\cite[Prop.\ 2.4]{HKMS15}.)

A more involved investigation is required when we consider the property of \emph{relative
computable categoricity} from computable structure theory.  A countable structure $\A$ is
relatively computably categorical if, for every two copies $\B$ and $\C$ of $\A$ with domain
$\omega$, there is an isomorphism from $\B$ onto $\C$ which is computable from the atomic
diagrams of $\B$ and $\C$.  An equivalent condition, discovered in \cite{AKMS89,Chis}
and expanded in \cite{M17}, is for $\A$ to have a Scott family of (finitary) $\Sigma_1$ formulas
which is $e$-reducible to the (finitary) $\Sigma_1$-theory of an expansion of $\A$
by finitely many constants.  This in turn is equivalent
to a uniform version of categoricity. Useful sources for details include \cite{AK00}
and \cite{K86}, which we will assume here as background.

\begin{thm}
\label{thm:cat^*fields}
The class
$$ \set{[F]\in\FF_0^*/\!\cong}{F\text{~is relatively computably categorical}}$$
of isomorphism types of relatively computably categorical fields in $\FF_0^*$
(in the language with all root predicates $R_n$)
has measure $1$ and is comeager within $\FF^*_0/\!\cong$.  Indeed, there is a single Turing functional
$\Theta$ such that
$$\set{[F]\in\FF_0^*/\!\cong}{(\forall K\cong F)~\Theta^{F\oplus K}:F\to K
\text{~is an isomorphism}]}$$
is comeager of measure $1$ there.
\end{thm}
\begin{pf}
We prove the stronger statement, writing $f=\Theta^{F \oplus K}$
for simplicity.  (We have now come to identify a field $F_D$ with its atomic diagram $D$,
writing $F$ for both and trusting the reader to understand which is intended.)
Let $f_0$ be the empty map.  On input $x\in F$, we
compute $f_{x+1}(x)$ using recursion on $x$, knowing $f_x(0),\ldots,f_x(x-1)$.
From the atomic diagram $F$, we can find
the minimal polynomial $p(0,\ldots,x-1,X)$ of each $x$ in $F$ over 
the subfield $F_x$ generated by $\dom{f_x}=\{ 0,\ldots,x-1\}$, 
and then find all roots $x_1=x,x_2,\ldots,x_k$ of $p(0,\ldots,x-1,X)$ in $F$.
(The root predicates compute $k$, so we know when we have them all.)
Likewise, we can then find all $k$ roots $y_1,\ldots,y_k$ of $p(f_x(0),\ldots,f_x(x-1),Y)$
in $K$; there must be exactly $k$ of them, since $F\cong K$ (and since, by induction,
$f_x=f\res x$ extends to some isomorphism).  If $k=1$, we define
$f_{x+1}(x)=y_1$ immediately.  Otherwise, for each $i\leq k$, we use the root predicates
in $F$ and $K$ to search for a polynomial
$q\in\Q[X_0,\ldots,X_x,Z]$ such that
$$(\exists a\in F)~q(0,\ldots,x,a)=0 \iff (\forall b\in K)~q(f_x(0),\ldots,f_x(x-1),y_i,b)\neq 0.$$
Whenever we find such a $q$, we say that it has \emph{ruled out $y_i$} as an image for $x$.
If we reach a stage where all but one of $y_1,\ldots,y_k$ has been ruled out thus, then
we define $f_{x+1}(x)$ to be the remaining $y_i$.  If this never happens, then the program
never halts.

It is clear that, if $f$ actually is total, then it will be an isomorphism from $F$ onto $K$.
Of course, there will be cases where $f$ is not total.  For example, whenever $F$ is a
proper normal extension of $\Q$, this $f$ will not be total.  (This includes many of the best-known
fields in $\FF_0$, such as $\Qbar$ and all splitting fields properly extending $\Q$.)
Also, if $f=\Theta^{F \oplus K}$ is non-total for some particular $F$ and $K$,
then $\Theta^{\Ftilde\oplus\Ktilde}$ will fail to be total whenever $F\cong\Ftilde$
and $K\cong\Ktilde$, since the image of $x$ in $\Ftilde$ will always be a stumbling block.

However, in order for $f(x)$ to diverge, two of the roots in $K$ -- say $y_1$ and $y_2$, without loss
of generality -- must both fail ever to be ruled out by the procedure above.  Now, for every
finite initial segment $\sigma\in 2^{<\omega}$ of $F$ and every conjugate $x_2$ of $x_1=x$
in $F$ over $\Q(0,\ldots,x-1)$, there is always a way to extend $\sigma$
so that some $q$ does distinguish between the two of them.
It follows that the class of fields $F$ for which (for some $K$, equivalently for all $K\cong F$)
two such $y_1$ and $y_2$ exist for this $x$ is a nowhere dense class.  Therefore, the class
of isomorphism types in $\FF_0$ for which $\Theta$ fails to work is a meager class.
Moreover, for any given $x=x_1$ and for each of its conjugates $x_i\neq x$ in $F$, the probability
that no $q$ ever distinguishes between them is $0$, since Lemma \ref{lemma:distinguish} will
give an infinite collection of polynomials, each of which has a one-half chance to distinguish
between them, and for which these one-half chances are all independent.
This completes the proof of the theorem.
\qed\end{pf}

In the foregoing proof, Lemma \ref{lemma:distinguish} establishes the intuitively
clear fact that, with probability $1$, two conjugates $y_1$ and $y_2$ in an algebraic field $F$
can always be distinguished by polynomials of the form $f(y_1,Z)$ and $f(y_2,Z)$:
for some $f$, one of these will have a root in the field and the other will not.
We will next prove the more surprising Theorem \ref{thm:catfields},
stating that, even without the root relations in the language,
measure-$1$-many of the fields in $\FF_0$ remain relatively computably categorical.
This will require Lemma \ref{lemma:distinguish} in exact detail.

\begin{lemma}
\label{lemma:distinguish}
Let $\alpha_0,\ldots,\alpha_n\in\Qbar$ be algebraic numbers conjugate over $\Q$.
Then, for every finite algebraic field extension $E\supseteq\Q$ with $\alpha,\beta\in E$,
and for all distinct $i\leq n$ and $j\leq n$, there exists an infinite set $D=\{ q_0,q_1,\ldots\}\subseteq\Q$
of rational numbers such that for every $k$, 
both of the following hold:
$$\sqrt{\alpha_i+q_k}\notin E(\sqrt{\alpha_i+q_l}, \sqrt{\alpha_j+q_l}~:~l\neq k)(\sqrt{\alpha_j+q_k});$$
$$\sqrt{\alpha_j+q_k}\notin E(\sqrt{\alpha_i+q_l}, \sqrt{\alpha_j+q_l}~:~l\neq k)(\sqrt{\alpha_i+q_k}).$$
Moreover, there is a procedure, uniform in $\alpha_0,\ldots,\alpha_n$, $i$, $j$, and the generators of $E$,
for deciding such a set $D$.
\end{lemma}
That is, adjoining $\sqrt{\alpha_i+q_k}$ to the field generated by all the other square roots
(from other elements $q_l\in D$, with $l\neq k$)
will not cause $\sqrt{\alpha_j+q_k}$ to appear in that same field, nor vice versa.  The same
therefore holds for adjoining $\sqrt{\alpha_i+q_k}$ to any subfield of that field.
\begin{pf}
We start by setting $E_0$ to be the normal closure of $E$,
the smallest extension of $E$ (clearly also of finite degree) which is normal over $\Q$.
We first use $E_0$ to find $q_0$, then define $E_1$ to be the normal closure
of $E_0(\sqrt{\alpha_i+q_0},\sqrt{\alpha_j+q_0})$ and repeat the process with
$E_1$ to get $q_1$, and so on, using the following step recursively.

We now appeal to the following well-known
fact, which appears, for instance, as \cite[Lemma 11.6]{FJ86}.
\begin{lemma}
\label{lemma:FJ}
Let $L$ be a finite separable extension of a field $K$, and let $f\in L(Z)[T]$ be irreducible
as a polynomial in $T$ over $L(Z)$.  Then there exists an irreducible $p\in K(Z)[T]$
such that for every $q\in K$, if $p(q,T)$ is irreducible in $K[T]$, then $f(q,T)$
is irreducible in $L[T]$.
\end{lemma}
We apply this lemma with $E_k$ as $L$ and with $\Q$ as $K$.  Recall that
$\Q$ is a \emph{Hilbertian field}:  for every $p\in\Q(Z)[T]$ irreducible over $\Q(Z)$,
there exist infinitely many $q\in\Q$ for which $p(q,T)$ is irreducible in $\Q[T]$.
Therefore, the lemma implies that for every irreducible $f\in E_k(Z)[T]$,
some $q\in\Q$ makes $f(q,T)$ irreducible in $E_k[T]$.  Now the polynomial
$T^2-\frac{\alpha_i+Z}{\alpha_j+Z}$ is irreducible in the polynomial ring $E_k(Z)[T]$ provided that
$\frac{\alpha_i+Z}{\alpha_j+Z}$ is not a square in $E_k(Z)$, which holds since $\alpha_i\neq\alpha_j$.
Therefore, there exists some $q\in\Q$ for which $T^2-\frac{\alpha_i+q}{\alpha_j+q}$
is irreducible in $E_k[T]$.  In particular, $(\alpha_i+q)$ and $(\alpha_j+q)$ are not
both squares in $E_k$.  However, if $(\alpha_i+q)$ were a square, say $\alpha_i+q=z^2$,
then we could apply an automorphism $h$ of $\Qbar$ mapping $\alpha_i$
to its $\Q$-conjugate $\alpha_j$.  Since $E_k$ is normal over $\Q$, it would
also contain $h(z)$, so $\alpha_j+q=h(\alpha_i+q)=(h(z))^2$ would also have been a square
in $E_k$, which is impossible.  Therefore $(\alpha_i+q)$ is not a square in $E_k$,
nor is $(\alpha_j+q)$, by the same reasoning.

Now suppose $\sqrt{\alpha_j+q}\in E_k(\sqrt{\alpha_i+q})$.  Then we have $x,y\in E_k$
with $(x+y\sqrt{\alpha_i+q})^2=\alpha_j+q$.  However, then
$$ x^2+2xy\sqrt{\alpha_i+q}+y^2(\alpha_i+q) = \alpha_j+q\in E_k,$$
so the coefficient $2xy$ must be zero.  But $y=0$ would force $\alpha_j+q$ to be a square in $E_k$,
while $x=0$ would force $\frac{\alpha_i+q}{\alpha_j+q}=\frac1{y^2}$ to be a square
in $E_k$.  Since both of these are impossible, we see that
$\sqrt{\alpha_j+q}\notin E_k(\sqrt{\alpha_i+q})$ and likewise $\sqrt{\alpha_i+q}\notin E_k(\sqrt{\alpha_j+q})$.

This argument shows that infinitely many numbers $q\in\Q$ exist for which
$\sqrt{\alpha_j+q}\notin E_k(\sqrt{\alpha_i+q})$ and $\sqrt{\alpha_i+q}\notin E_k(\sqrt{\alpha_j+q})$.
Knowing the generators of $E_k$, we have a splitting algorithm for $E_k$,
by Kronecker's Theorem, and so we may identify such a $q$ when we find one.
Choose $q_k$ to be some such $q$ which (as an element of the domain $\omega$ of $\Qbar$)
is greater than $k$; this will make our set $D$ decidable.
Setting $E_{k+1}$ to be the normal closure of
$E_k(\sqrt{\alpha_i+q_k},\sqrt{\alpha_j+q_k})$, we proceed by recursion on $k$,
and thus build the set $D$.

Now suppose for a contradiction that for some $k$,
$$\sqrt{\alpha_j+q_k}\in E(\sqrt{\alpha_i+q_l}, \sqrt{\alpha_j+q_l}~:~l\neq k)(\sqrt{\alpha_i+q_k}).$$
Take any $k$ for which this holds, and fix the least $p$ such that
$$\sqrt{\alpha_j+q_k}\in E(\sqrt{\alpha_i+q_l}, \sqrt{\alpha_j+q_l}~:~l\leq p~\&~l\neq k)(\sqrt{\alpha_i+q_k}).$$
By our construction, $\sqrt{\alpha_j+q_k}\notin E_k(\sqrt{\alpha_i+q_k})$,
so clearly $p>k$, and (by the minimality of $p$) either the adjoinment of $\sqrt{\alpha_i+q_p}$ to the field
$$ \Etilde = E(\sqrt{\alpha_i+q_l}, \sqrt{\alpha_j+q_l}~:~l < p~\&~l\neq k)(\sqrt{\alpha_i+q_k}) $$
or the subsequent adjoinment of $\sqrt{\alpha_j+q_p}$ to $\Etilde(\sqrt{\alpha_i+q_p})$
caused $\sqrt{\alpha_j+q_k}$ to enter the field.
But each of these two extensions was of degree $2$, and so that extension
must also be generated by $\sqrt{\alpha_j+q_k}$.  Therefore either
$$ \sqrt{\alpha_i+q_p}\in\Etilde(\sqrt{\alpha_j+q_k})\subseteq E_p$$
or
$$\sqrt{\alpha_j+q_p}\in \Etilde(\sqrt{\alpha_i+q_p},~\sqrt{\alpha_j+q_k})\subseteq E_p(\sqrt{\alpha_i+q_p}),$$
both of which contradict our construction at stage $p$.
This proves one of the two conditions required by the lemma, and
by symmetry on $i$ and $j$, the other also holds.
\qed\end{pf}

\begin{thm}
\label{thm:catfields}
The class
$$ \set{[F^*]\in\FF^*_0/\!\cong}{F\text{~is relatively computably categorical in }\FF_0}$$
of isomorphism types of relatively computably categorical fields in $\FF_0$
(in the language without the root relations $R_n$)
has measure $1$ and is comeager within $\FF_0^*/\!\cong$.  
\end{thm}
To be clear:  the measure here is still taken in the isomorphism space
$\FF_0^*/\!\cong$ with the root relations, as this is the only isomorphism
space in which this measure makes sense.  The theorem states that measure-$1$-many
of these fields remain relatively computably categorical even when the
isomorphisms are forbidden to use information about the root relations
from the atomic diagrams.  However, we lose the uniformity of Theorem \ref{thm:cat^*fields}.
\begin{pf}
Fix any $\ep>0$.  We will enumerate, effectively, a family $\Sigma$ of formulas
which, for all $[F]$ in a subset of $\FF_0^*/\!\cong$ of measure $>1-\ep$, forms
a Scott family for $F$.  This will show that all these $[F]$ are relatively computably
categorical, since having a c.e.\ Scott family implies relative computable categoricity.
Perhaps the most striking thing about this family is that it will be built with no reference
to any specific $[F]$ at all, apart from its use of a fixed computable presentation of $\Qbar$,
which we denote simply as $\Qbar$, and a fixed computable enumeration $p_0,p_1,\ldots$
of the irreducible monic polynomials in $\Q[X]$, with the coefficients viewed as elements
of $\Qbar$.

Consider the $n$-th polynomial $p_n(X)$.  Let $d$ be its degree.  If $d=1$,
so that $p_i=X-\frac{m}{k}$ for some $m\in\Z$ and nonzero $k\in\omega$,
then we add to $\Sigma$ a formula saying that $kX=m$.
Otherwise, we now apply Lemma \ref{lemma:distinguish} with a vengeance,
with the splitting field of $p_n$ over $\Q$ as our $E$, and with the roots
$\alpha_1,\ldots,\alpha_d$ of $p_n$ in $\Qbar$.  Let $\del=\frac{\ep}{2^{n+1}}$,
and compute an $N$ so large that both of the following hold:
\begin{itemize}
\item
The probability that, out of $100 N$ independent coin flips, at least $40N$ will be heads,
is $>(1-\del)^{\frac1{2d}}$; and
\item
The probability that, out of $100N$ independent flips of two coins each,
at most $35N$ will yield two heads, is $>(1-\del)^{\frac1{d(d-1)}}$.
\end{itemize}
Fix elements $q_0,\ldots,q_{100N-1}$ from the set $D$ given by
Lemma \ref{lemma:distinguish}.  For each subset $S\subseteq
\{ q_0,\ldots,q_{100N-1}\}$ with $|S|=40N$, we add to $\Sigma$ the fomula
$$ p_n(X)=0~\&~\bigwedge_{q\in S}~\exists Y_q~Y^2=X+q,$$
saying that $X$ is one of the $\alpha_i$ and that, for all $q\in S$,
the field contains an element $\sqrt{\alpha_i+q}$.

The first item ensures that, with probability $>\sqrt{1-\del}$, every $\alpha_i$ in $F$
will satisfy $(\exists x)~x^2=\alpha_i+q_k$ for at least $(40N)$-many values $j<100N$.
Therefore, with this probability, every root of $p_n$ in $F$ (if there are any) will indeed
realize one of the formulas now in $\Sigma$.  On the other hand, the second item ensures
that the probability is $>((1-\del)^{\frac1{d(d-1)}})^{\binom d2}=\sqrt{1-\del}$ that, for every two distinct roots $\alpha_i,\alpha_j$
of $p_n$ in $F$, no more than $(35N)$ values of $k<100N$ have the property that
$$ (\exists y\in F)~y^2=\alpha_i+q_k~~~\&~~~(\exists z\in F)~z^2=\alpha_j+q_k.$$
Therefore, with probability $>1-\del = 1-\frac{\ep}{2^{n+1}}$, $\Sigma$ behaves exactly
as a Scott family should, as far as the roots of $p_n$ are concerned.  (No other formula
in $\Sigma$ can be realized by any root of $p_n$, since all such formulas either specify
a rational value for $X$, or else require $X$ to be a root of some $p_m$ with $m\neq n$.)

Finally, to make a true Scott family, we adjoin to $\Sigma$ all conjunctions of the form
$$ \bigwedge_{i\leq m}~\alpha_{n_i}(X_{n_i}),$$
for all strictly increasing tuples $(n_0,\ldots,n_m)\in\omega^{<\omega}$,
where $\alpha_{n_i}$ can be any of the formulas added to $\Sigma$ for $p_{n_i}$ above.
Scott families, after all, must contain formulas realizing every $n$-tuple from a structure,
even though, in this case, the single-variable formulas actually suffice, being specific
enough that, when required to identify a new element with respect to the finitely many
already considered, they can actually identify it uniquely without even considering
the preceding elements.

The probability that our $\Sigma$ fails to work for the roots of a single $p_n$
is $<\frac{\ep}{2^{n+1}}$, by the analysis above.  Therefore, for at least $(1-\ep)$-many
of the isomorphism types in $\FF^*_0/\!\cong$, $\Sigma$ is a computably enumerable
Scott family, and thus all these types are relatively computably categorical.  Moreover,
$\Sigma$ does not use the root predicates anywhere, and therefore these isomorphism
types remain relatively computably categorical even in the original language of fields,
without the root predicates. Since $\ep$ was arbitrary, the theorem follows.
\qed\end{pf}

It is well-known that normal algebraic field extensions of $\Q$ are also relatively
computably categorical, even without the root predicates, although the specific
procedure is different from the functionals in Theorems \ref{thm:cat^*fields} and
\ref{thm:catfields}.  With normal extensions, one simply finds the minimal polynomial
$p$ of $x$ over $\Q(0,1,\ldots,x-1)$ and defines $f_{x+1}(x)$ to be the first root of
$p(f_x(0),\ldots,f_x(x-1),Y)$ that appears in the target field $K$, knowing
that by normality, this map must extend to an isomorphism.  Of course, normal
extensions are a meager class of measure $0$ within the class of all algebraic
extensions of $\Q$; our point is simply that separate procedures may work
in special cases.  In the forthcoming work \cite{FM19}, it is shown that the procedures
given in Theorems \ref{thm:cat^*fields} and \ref{thm:catfields} succeed for all
\emph{random fields}, i.e., for all $F$ such that the $h\in 2^\omega$ corresponding
to $[F]$ is Martin-L\"of random, or even just Schnorr random.

\comment{
There are situations in which one can do a similar analysis not on isomorphism types,
but on classes of countable structures.  For example, the class $\LL$ of all (infinite)
countable linear orders is also homeomorphic to Cantor space:  one visualizes the linear order
$<$ being built by deciding first whether $0<1$, next whether $2<0$, next whether $2<1$ (unless
this was already determined by the preceding two answers -- which is decidable, of course),
next whether $3<0$, and so on.  In this particular example, however, the class $\LL$ is dominated
by a single isomorphism type:  it is readily seen that measure-$1$ many (and comeager-many)
orders in $\LL$ are isomorphic to the countable dense linear order without end points.  It follows that,
in $\LL$ (\emph{not} modulo isomorphism), the analogue of Theorem \ref{thm:catfields} holds,
whereas in $\LL/\!\cong$ it surely fails, under any reasonable notion of measure or category,
since only countably many isomorphism types of linear orders are relatively computably
categorical at all.
}

\section{Different Measures}
\label{sec:Haar}

In the preceding section, we transferred the Lebesgue measure on Cantor space
to the space $\FF_0^*/\!\cong$, using our homeomorphism from Theorem \ref{thm:algfields}.
This seems like a natural choice, yet in this section we present another measure on
$\FF_0^*/\!\cong$, which we consider to be, if anything, a more natural choice.
The results of Section \ref{sec:measure} hold no matter which of the two measures one uses.

We were careful, in the proof of Theorem \ref{thm:algfields}, to consider only polynomials
$f_n$ of prime degree.  This avoided one possible trap regarding the measure:  it does not
depend on the particular enumeration $f_0, f_1,\ldots$ of the monic polynomials in $E[X]$ there.
Had we not required prime degree, the following situation could have arisen.
\begin{enumerate}
\item[(a)]
If $f_0(X)=X^2-2$ , then the measure of the class of all isomorphism types of fields
containing a square root of $2$ is $\frac12$.
\item[(b)]
However, if instead $f_0(X)=X^4-2$ and $f_1(X)=X^2-2$, then Lebesgue measure would
dictate that half of all fields contain a fourth root of $2$, and that half of the remaining
fields contain a square root of $2$ (but no fourth root).  In this case, the measure of the
class of all isomorphism types of fields containing a square root of $2$ is $\frac34$.
\end{enumerate}
By requiring prime degree, we avoided this trap to some extent.
Nevertheless, certain peculiarities
of the measure remain.  For example, the fields containing an $181$-st root of $2$ form a class of
measure $\frac12$, since $181$ is prime, whereas the fields containing
a $180$-th root of $2$ form a class of far smaller measure, since such a root
is built up by a long series of prime-degree extensions.  This is not inconsistent in any way,
but it seems a little strange, and there is a way of avoiding it.

\begin{defn}
\label{defn:Haar}
The \emph{Haar-compatible measure} $H$ of a basic open set $\U_\sigma$ in Cantor space
is defined as follows, using the same notation as in the proof of Theorem \ref{thm:algfields}:
within a computable copy $E$ of $\Qbar$, we choose the same polynomials $f_\sigma$
and build the same subfields $F_\sigma$ of $E$.  For the empty string $\lambda$,
of course we define $H(\U_\lambda)=1$.  Given $H(\U_\sigma)$, we set
$$ H(\U_{\sigma\widehat{~}1}) = \frac1d\cdot H(\U_\sigma)~~~~\&~~~~
H(\U_{\sigma\widehat{~}0}) = \frac{d-1}d\cdot H(\U_\sigma),$$
where $d$ is the degree of the polynomial $f_\sigma$.  The Haar-compatible measure
of an arbitrary subset $\V\subseteq 2^\omega$ is then the infimum
of the sums $\Sigma_i H(\U_{\sigma_i})$, over all countable sequences
$\la \sigma_i\ra_{i\in\omega}$ for which $\V\subseteq\bigcup_i\U_{\sigma_i}$.
\end{defn}
So the set $\U$ of fields containing an $11$-th root of $2$ has 
$H(\U)=\frac1{11}$, as opposed to its Lebesgue measure $\frac12$.
Of course, $H(\U_{\sigma\widehat{~}0})+H(\U_{\sigma\widehat{~}1})=H(\U_\sigma)$,
so this really is a measure.  The next lemma justifies the name ``Haar-compatible measure.''

\begin{lemma}
\label{lemma:Haar}
For every finite Galois extension $K$ of $\Q$, the set $\V$
of all fields in $\FF_0^*/\!\cong$ containing $K$ has
$H(\V)$ equal to the usual Haar measure $\frac1{[K:\Q]}$
of the (pointwise) stabilizer of $K$ within the Galois group $\Gal{\Qbar}{\Q}$.
\end{lemma}
\begin{pf}
Let $h$ be any element of Cantor space, with corresponding field
$F_h = \cup_{\sigma\subset h}F_\sigma\in\FF_0^*$.
Assuming $K\neq\Q$, there will be at least one $\sigma\subset h$, for which $F_\sigma\cap K$
is a proper subfield of $F_{\sigma\widehat{~}1}\cap K$.  These $\sigma$ are the nodes
at which $h$ makes decisions about whether to include $K$ (if $h(|\sigma|)=1$)
or not (if $h(|\sigma|)=0$).  At the first such $\sigma_1$,
with $p_1=\deg{f_{\sigma_1}}$, we will have $H(\U_{\sigma_1\widehat{~}0})=
\frac{p_1-1}{p_1}\cdot H(\U_{\sigma_1})$, and, since $K$ is normal over $\Q$,
none of these fields will lie in $\V$.  (Notice that if $K$ were not normal,
it might still be possible to have an $h\supset\sigma_1$ such that $F_h$ contained
a subfield isomorphic to $K$.)
For the remaining $\frac1{p_1}\cdot H(\U_{\sigma_1})=H(\U_{\sigma_1\widehat{~}1})$-many
fields, either $F_{\sigma_1\widehat{~}1}\supseteq K$ and these fields all lie in $\V$;
or else we continue up through Cantor space until we find another $\sigma_2$
with the same property.  After we reach $\sigma_n$ (where $n$ is the number of prime factors
of $[K:\Q]$, counted by multiplicity), we must have $K\subseteq F_{\sigma_n\widehat{~}1}$,
and multiplying all the numbers to this point proves the lemma, since the degrees
of all the extensions $[F_{\sigma_i\widehat{~}1}:F_{\sigma_i}]$ must have product
$[K:\Q]$.
\qed\end{pf}

Lemma \ref{lemma:Haar} emphatically does require that $K$ be normal over $\Q$.
Our Haar-compatible measure for finite non-normal field extensions is more complicated,
as is the Lebesgue measure.  For example, consider the class $\V$ of fields
$K$ containing a fourth root of $2$ (equivalently, a copy of $\Q(\sqrt[4]2)$).
If our sequence of polynomials in the proof of Theorem \ref{thm:algfields}
begins with:
$$f_0 = X^2-2;~~~~~f_1=X^2-\sqrt2;~~~~~f_2=X^2+\sqrt2,$$
then $H(\V)=\frac38$, as $\V$ corresponds to the union of $\U_{11}$
and $\U_{101}$.  However, if instead we used a sequence beginning with:
$$g_0=X^2+1;~~~~~g_1=f_0;~~~~~g_2=f_1;~~~~~g_3=f_2,$$
then $H(\V)= \frac5{16}$ instead, as now $\V$ corresponds to the union
$\U_{111}\cup\U_{011}\cup\U_{0101}$.  (The point is that the string
$\sigma=110$ rules out any square root of $-\sqrt2$, since the presence
of $i$ in the field means that we must have either no fourth roots of $2$,
or else all four of them.)  The Haar-compatible measure is not to blame, insofar as
$H(\V)$ is also the Lebesgue measure of $\V$ in both of these cases.
For non-normal field extensions, the measure still depends on the
sequence of polynomials used.  We conjecture that various questions
arising in this and the following section are more readily addressed
if one considers only the class of normal algebraic field extensions
of $\Q$, in place of the class $\FF_0$.

\section{Topologies on $2^\omega$}
\label{sec:Scott topology}

As promised above, we return here to the topology of $\FF_0/\!\cong$,
which arose in Section \ref{sec:fields}.  There we noted that it is not homeomorphic to
Cantor space.  (By \emph{Cantor space} we really mean not just the set $2^\omega$,
but rather this set with its usual topology, where the basic open sets are defined by
$\U_\sigma=\set{h\in 2^\omega}{\sigma\subset h}$ for all $\sigma\in 2^{<\omega}$.)
Here we consider other possibilities.

The \emph{Scott topology} on $2^\omega$ is a sort of positive version of Cantor space.
Here the basic open sets are defined to be the sets
$$\W_F = \set{h\in 2^\omega}{F\subseteq h^{-1}(1)}$$
for all finite subsets $F$ of $\omega$.  If $F=\{ 3,10\}$, then $\W_F$
contains all $h\in 2^\omega$ with $h(3)=h(10)=1$, or equivalently,
all subsets of $\omega$ which contain both $3$ and $10$.  However,
sets such as $\set{h\in 2^\omega}{h(3)=0}$, which were open in
Cantor space, are not open in the Scott topology.  We will write
$\SS$ for the space $2^\omega$ under the Scott topology.

The Scott topology seems much closer in nature to $\FF_0/\!\cong$, being
defined by positive information only.  Both contain an element which
lies in every nonempty open set (for $\SS$, it is the constant function $h=1$)
and an element which sits in no proper open subset of the space (the other
constant function $h=0$).  In fact, there is a very natural equivalence relation
on $2^\omega$ whose quotient has the Scott topology.  Earlier we used $=^e$
to denote a binary relation on $\omega$.  Now we revise it to denote the related
notion on $2^\omega$:
$$ A =^e B \iff (\forall x)[(\exists y~\la x,y\ra\in A) \iff (\exists z \la x,z\ra\in B)].$$
This relation is best understood by considering the image $\pi_1(A)=\set{x\in\omega}{\exists y~\la x,y\ra\in A}$
of $A$ under projection $\pi_1$ onto the first coordinate.  The set
$\pi_1(A)$ is said to be \emph{enumerated by $A$}, and $A$ and $B$ satisfy $=^e$
if and only if they enumerate the same set.  Notice that a set is $C$-computably
enumerable just if there is a $C$-computable set enumerating it:  the terminology captures
our intuition about enumeration of sets.  Likewise, it gives an equivalent definition
of enumeration reducibility $\leq_e$: an \emph{enumeration reduction
from $D$ to $C$} is simply a Turing functional $\Gamma$ such that, for every
enumeration $B$ of $D$, $\Gamma^B$ is the characteristic function of an enumeration
of $C$.

The homeomorphism from $2^\omega/\!=^e$ onto $\SS$ is simple: one maps each $A\in 2^\omega$
to $\pi_1(A)$ in $\SS$.  The preimage of each $\W_F$ is open, so this is continuous,
and it respects $=^e$, so it gives rise to a
continuous function from $2^\omega/\!=^e$ onto $\SS$, which is quickly seen to be
bijective.  The inverse is even easier: map $C\in\SS$ to the $=^e$-class of
the set $\set{\la x,0\ra}{x\in C}$.  Again, this is easily seen to be continuous,
so we have a homeomorphism.  Indeed, the natural bases for the two spaces both
consist of sets defined by the requirement of including a certain finite subset $F$ of $\omega$.

We now introduce another natural and closely related equivalence relation $=^f$ on $2^\omega$,
which aids in effective classification of several classes of countable structures
not considered in this article, such as equivalence structures and torsion-free abelian groups
of rank $1$.  Recall that the $x$-th \emph{column} of $A\subseteq\omega$ is defined
to be $A^{[x]}=\set{y\in\omega}{\la x,y\ra\in A}$.  The relation $=^f$ is defined by:
\begin{align*}
A =^f B &\iff (\forall x)[A^{[x]}~\Ecard~B^{[x]}]\\
&\iff (\forall x)~|\set{y}{\la x,y\ra\in A}| = |\set{y}{\la x,y\ra\in B}|.
\end{align*}
So we think of each column $A^{[x]}$ of $A$ as ``counting'' something -- perhaps the number
of times a fixed element of a torsion-free abelian group is divisible by the prime $p_x$, for example.
Instead of enumerating a set, $A$ thus approximates a function $g_A$, by computable approximations
from below.  Now $A=^f B$ just if $g_A=g_B$
as functions.  Notice, however, that this is not quite the same notion as the limitwise monotonic
approximation of functions used in many other topics in computable structure theory, since the
function $g_A$ maps $\omega$ into $(\omega+1)$, taking on
the value $\omega$ whenever the set $A^{[x]}$ is infinite.

Basic open sets in $2^\omega/\!=^f$ may be thought of as given by finite pieces of functions.
For instance, using the function $g_0$ with $g_0(3)=4$, $g_0(10)=2$, and $g_0(x)=0$ for all
other $x$, we get a basic open set containing all $A$ for which $g_A(3)\geq 4$ and $g_A(10)\geq 2$.
Notice that $g_0$ must have finite support (i.e., $g_0^{-1}(0)$ must be cofinite) and must take on
only finite values:  setting $g_0(x)=\omega$ is not allowed.  (A finite initial segment of $A$
can ensure $g_A(3)\geq 4$, but no finite initial segment can ensure $g_A(x)=\omega$.)
The space $2^\omega/\!=^f$ is similar in many respects to the Scott topology, but the two
are not the same.

\begin{thm}
\label{thm:Scott}
The spaces $2^\omega/\!=^e$ (i.e., the Scott topology) and $2^\omega/\!=^f$ are not homeomorphic.
\end{thm}
\begin{pf}
We use the notion of a \emph{principal open set}.  In a topological space, an open set $\V$
is \emph{principal} if there exists some point $x$ in the space for which $\V$ is the smallest
open set containing $x$.  (Such an $x$ may be said to \emph{generate} this $\V$.)
This property is clearly preserved under homeomorphisms.

Now in $2^\omega/\!=^e$ and $2^\omega/\!=^f$, principal open sets do exist, and in fact they
are just the basic open sets described above.  Intuitively, each one is generated by a given finite
subset of $\omega$.  Formally, for each finite $F\subset\omega$, the open set
$$\V_F = \set{[A]\in 2^\omega/\!=^e}{F\subseteq\pi_1(A)}$$
is principal, being generated by $[F]$, and all principal open sets are of this form, including
the entire space, which is generated by the singleton $=^e$-class $[\emptyset]$.
(In the notation for the Scott topology $\SS$, this $\V_F$ was the set $\W_F$.)

Likewise, in $2^\omega/\!=^f$, principal open sets are of the form
$$\X_g = \set{[A]}{(\forall x) g_A(x)\geq g(x)},$$
where $g$ ranges over all total functions from $\omega$ into $\omega$
with finite support.  $\X_g$ is generated by the $=^f$-class
$[\set{\la x,y\ra}{y<g(x)}]$ of a set whose $x$-th column always has cardinality $g(x)$.

Now we can explain why the two spaces are not homeomorphic.  Let $\V=\V_F$
be any principal open set in $2^\omega/\!=^e$.  Then the principal open supersets
of $\V$ in this space are precisely those of the form $\V_G$ with $G\subseteq F$,
and so there are exactly $2^{|F|}$ principal open sets containing $\V_F$,
including both $\V_F$ itself and the entire space.  However, in $2^\omega/\!=^f$,
the principal open supersets of each $\X_g$ are those $\X_h$ satisfying
$(\forall x~h(x)\leq g(x))$.  In particular, the principal open set $\X_{g_2}$
for the function $g_2(0)=2$ (with support $\{2\}$)
has exactly three principal open supersets:  itself, the entire space, and
$\X_{g_1}$, where $g_1(0)=1$.  Thus no homeomorphism
can map $\X_g$ to any principal open set in $2^\omega/\!=^e$.
\qed\end{pf}

Theorem \ref{thm:Scott} is of interest on its own, but another reason for proving it was to
prepare the way for proving the answer to our initial question from this section.
\begin{thm}
\label{thm:Alg}
The Scott topology and $\FF_0/\!\cong$ are not homeomorphic; nor are $2^\omega/\!=^f$ and $\FF_0/\!\cong$.
\end{thm}
\begin{pf}
We use the characterization from Theorem \ref{thm:Scott} of the principal open sets
in $2^\omega/\!=^e$ and $2^\omega/\!=^f$.  In $\FF_0/\!\cong$, by much the same reasoning,
the principal open sets are the basic open sets from Lemma \ref{lemma:FF_0}.
As an example,
consider the subfield $\Q(\theta)$ of $\Qbar$ generated by a single primitive fifth root $\theta$
of unity.  The minimal polynomial of $\theta$ is the cyclotomic polynomial
$X^4+X^3+X^2+X+1$, and its conjugates are $\theta^2$, $\theta^3$, and $\theta^4$,
all of which lie in $\Q(\theta)$.  Thus the Galois group of this normal extension is $\Z/(4)$.
Since this Galois group has a proper nontrivial subgroup, there is a subfield of degree $2$
within $\Q(\theta)$, generated by an element $\sqrt{z}$ for some $z\in\Q$.  Now the subfield
$\Q(\sqrt{z})$ generates a principal open set (more formally, the $\cong$-class $[\Q(\sqrt{z})]$
generates it), which is a proper superset of the principal open set generated by $[\Q(\theta)]$.
$[\Q]$ itself generates another principal open superset, namely the entire space, and since
$\Q(\theta)$ has no other subfields, the principal open set generated by $[\Q(\theta)]$ has
exactly three principal open supersets.  As before, this shows that $2^\omega/\!=^e$
cannot be homeomorphic to $\FF_0/\!\cong$.

For $=^f$, a similar numerical argument will not work.  Instead, we consider the principal
open supersets of a given principal open set as a partial order, under $\subset$, and show
that certain finite partial orders possible in $\FF_0/\!\cong$ are impossible in $2^\omega/\!=^f$.
In the latter space, when $\X_g$ is a principal open set, its principal open supersets
are those of the form $\X_h$ with $h\leq g$ on all inputs.  In particular, let
$\omega-g^{-1}(0) = \{ x_1<\cdots<x_n\}$ and, for each $i<n$, let
$$ g_i(x)=\left\{\begin{array}{cl}g(x),&\text{~if~}x=x_i;\\0,&\text{~if not.}\end{array}\right.$$
Each $g_i$ gives rise to a proper chain $\X_{i,0}\supset\X_{i,1}\supset\cdots\supset\X_{i,g(x_i)}$
of principal open supersets of $\X_g$, with $\X_{i,0}$ always being the entire space.
The principal open supersets of $\X_g$ are just the $\Pi_i (1+g(x_i))$-many
intersections of these finitely many sets, and those which are maximal
among the principal open proper subsets of $2^{\omega}/\!\cong$ are
precisely the sets $\X_{i,1}$, for $i=1,\ldots,n$.

Now we wish to produce a field $K$, algebraic of finite degree over $K$, for which the
partial order of the subfields of $K$ is not of the form above.  Our specific example is the splitting
field $K$ of the polynomial $X^4-2$, for which we present here the lattice of subfields
(up to isomorphism) under inclusion:

\hspace{2cm}
\setlength{\unitlength}{0.25in}
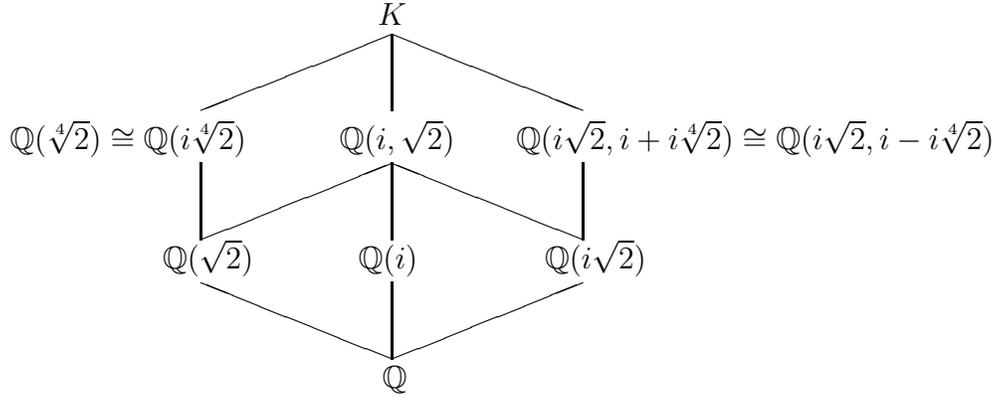
\begin{picture}(11,9)(1.5,0.2)

\put(5.3,0.6){$\Q$}

\put(0.7,3.1){$\Q(\sqrt2)$}
\put(4.8,3.1){$\Q(i)$}
\put(8.7,3.1){$\Q(i\sqrt2)$}


\put(1.5,3.7){\line(0,1){1.6}}
\put(5.5,3.7){\line(0,1){1.6}}
\put(9.5,3.7){\line(0,1){1.6}}

\put(1.5,3.7){\line(5,2){4}}
\put(9.5,3.7){\line(-5,2){4}}

\put(-2.5,5.6){$\Q(\sqrt[4]2)\cong\Q(i\sqrt[4]2)$}
\put(4.4,5.6){$\Q(i,\sqrt2)$}
\put(8.1,5.6){$\Q(i\sqrt2,i+i\sqrt[4]2)\cong\Q(i\sqrt2,i-i\sqrt[4]2)$}

\put(5.2,8.2){$K$}

\put(5.5,1.2){\line(0,1){1.6}}
\put(5.5,1.2){\line(-5,2){4}}
\put(5.5,1.2){\line(5,2){4}}

\put(5.5,8){\line(0,-1){1.6}}
\put(5.5,8){\line(-5,-2){4}}
\put(5.5,8){\line(5,-2){4}}

\end{picture}

\noindent
One sees from this lattice that the principal open set $\U_{[K]}$
generated by $[K]$ in $\FF_0/\!\cong$ has eight
principal open supersets, including itself and the entire space.  Those generated
by $[\Q(\sqrt2)]$, $[\Q(i)]$, and $[\Q(i\sqrt2)]$ are maximal proper.  In order for an $\X_g$ as above
to give rise to the same lattice, $g$ would need to have exactly three elements
in its support (that is, $n=3$), and to have $(1+g(x_1))\cdot (1+g(x_2))\cdot (1+g(x_3))=8$.
Moreover, in order for the lattice not to have a chain of length $>3$, we would need
all $g(x_i)\leq 2$, making $g(x_i)=1$ for each $i$.
However, in the lattice for such an $\X_g$, each maximal proper subset
contains itself, $\X_g$, and exactly two of the other principal open sets:
for example, if $h(x_1)=1$ and $h=0$ everywhere else, then $\X_h$
contains $\X_f$ for exactly those $f$ with $h\leq f\leq g$, and such an $f$
must have $f(x_1)=1$, $f(x_2)\in\{0,1\}$, $f(x_3)\in\{ 0,1\}$, and $f=0$ elsewhere.
Thus, such an $\X_g$ cannot be mapped onto $\U_{[K]}$ by any homeomorphism,
as no $\X_h$ could be mapped onto $\U_{[\Q(i)]}$.
This proves the theorem.
\qed\end{pf}

The broader conclusion of the foregoing proof is that the question
of the actual structure of the principal open sets of $\FF_0/\!\cong$
is the same question as the structure of the finite subfields of $\Qbar$ under inclusion,
and thus related to the study of $\Gal{\Qbar}{\Q}$, a notoriously difficult topic.  We consider
this to be a fair validation of the claim in Section \ref{sec:fields} that $\FF_0/\!\cong$
is not a recognizable topological space, and therefore should not be regarded as
a useful effective classification of $\FF_0$ up to isomorphism.  Adding the root predicates
to the language really is necessary.  On the other hand, we would be quite willing to
accept spaces such as $2^\omega/\!=^e$, or modulo $=^f$, as giving effective classifications
of other classes $\C$ of structures, provided that a homeomorphism (preferably computable)
exists from that space onto $\C/\!\cong$.  Indeed, descriptive set theory provides quite
a number of standard Borel equivalence relations on Cantor space and on Baire space,
noncomputable but of relatively low complexity; these include the relations $E_0$, $E_1$,
$E_2$, $E_3$, $\Eset$, and $Z_0$ (see \cite{MSEALS} or many other places for their
definitions), and also the relation $\Eperm$, computably bireducible with $\Eset$,
which holds of $A$ and $B$ just if the columns of $A$ are exactly the columns of $B$,
up to a permutation of columns.  ($\Eset$ merely requires that every column of each
should also appear as a column in the other, disregarding multiplicities.)  Beyond that,
it is possible to combine these notions:  for example, one might mix $=^e$ with $E_0$
by declaring $A$ and $B$ to be equivalent whenever $\pi_1(A)~E_0~\pi_1(B)$ holds.
Another variant, the equivalence relation $\Ecard^{\forall}$ defined by
$$A~\Ecard^{\forall}~B \iff (\omega-\pi_1(A))~\Ecard~(\omega-\pi_1(B)),$$
is quickly seen to have $2^\omega/\Ecard^{\forall}$ homeomorphic to $\AA_0/\!\cong$
(just in the language of fields, with no dependence relations).

Finally, in addition to Cantor space, Baire space and all their quotients, there are the
topological spaces $\mathbb{R}$ and $[0,1]$ under the usual topologies.  We have no idea whether
any isomorphism space $\C/\!\cong$ might be homeomorphic to either of these:  it seems
unlikely, but for now it stands as a good question.

\parbox{4.7in}{
{\sc
\noindent
Department of Mathematics \hfill \\
\hspace*{.1in}  Queens College -- C.U.N.Y. \hfill \\
\hspace*{.2in}  65-30 Kissena Blvd. \hfill \\
\hspace*{.3in}  Queens, New York  11367 U.S.A. \hfill \\
Ph.D. Programs in Mathematics \& Computer Science \hfill \\
\hspace*{.1in}  C.U.N.Y.\ Graduate Center\hfill \\
\hspace*{.2in}  365 Fifth Avenue \hfill \\
\hspace*{.3in}  New York, New York  10016 U.S.A. \hfill}\\
\medskip
\hspace*{.045in} {\it E-mail: }
\texttt{Russell.Miller\at {qc.cuny.edu} }\hfill \\
}

\end{document}